\documentclass{amsart}

\usepackage{amssymb}
\usepackage{amsthm}
\usepackage{cases}
\usepackage{amsfonts}
\usepackage{mathrsfs}
\usepackage{amsmath}
\usepackage{extarrows}
\usepackage{color}
\newtheorem{theorem}{Theorem}[section]
\newtheorem{corollary}[theorem]{Corollary}
\newtheorem{lemma}[theorem]{Lemma}
\newtheorem{proposition}[theorem]{Proposition}

\numberwithin{equation}{section}

\newcommand{\RePt}{\mathrm{Re}\,}
\newcommand{\ImPt}{\mathrm{Im}\,}
\newcommand{\ball}{\mathbb{B}}

\newcommand{\calU}{\mathcal{U}}

\newcommand{\calS}{\mathcal{S}}
\newcommand{\bfi}{\mathbf{i}}

\newcommand{\bfrho}{\boldsymbol{\rho}}
\newcommand{\calB}{\mathcal{B}}
\newcommand{\bbB}{\mathbb{B}}
\newcommand{\bbC}{\mathbb{C}}

\begin{document}

\title[BMO and Hankel operators]{BMO and Hankel operators on Bergman space of the Siegel upper half-space}

\thanks{This work was supported by the Hainan Provincial Natural Science Foundation of China grant 120QN177.}

\author{Jiajia Si}
\email{sijiajia@mail.ustc.edu.cn}

\address{School of Science, 
Hainan University, 
Haikou, Hainan 570228, 
People’s Republic of China.}

\subjclass[2010]{Primary 47B35; Secondary 32A36}

\begin{abstract}
On the setting of the Siegel upper half-space we study the spaces of bounded and vanishing mean oscillations which are defined in terms of the Berezin transform, and we use them to characterize bounded and compact Hankel operators on Bergman space.

\end{abstract}

\keywords{Hankel operator; Bergman space; BMO; Berezin transform; Siegel upper half-space.}

\maketitle

\section{Introduction}

Let $\mathbb{C}^n$ be the $n$-dimensional complex Euclidean space.
For any two points $z=(z_1,\cdots,z_n)$ and $w=(w_1,\cdots,w_n)$ in $\mathbb{C}^n$ we write
\[
z\cdot\overline{w} := z_1\overline{w}_1+\cdots+z_n\overline{w}_n,
\]
and
\[
|z|:= \sqrt{z\cdot\overline{z}} = \sqrt{|z_1|^2+\cdots+|z_n|^2}.
\]
The set
\[
\calU=\left\{ z\in\mathbb{C}^n:\ImPt z_n>|z^{\prime}|^2\right\}
\]
is the Siegel upper half-space. Here and throughout the paper, we use the notation
\[
z=(z^{\prime},z_n),\,\,\,\, \text{where}\, z^{\prime}=(z_1,\cdots,z_{n-1})\in\mathbb{C}^{n-1}\,\, \text{and}\,\, z_n\in\mathbb{C}^{1}.
\]
The Bergman space $A^2(\calU)$ is the space of all complex-valued holomorphic functions $f$ on $\calU$ such that 
\[
\|f\|:=\left(\int_{\calU} |f(z)|^2 dV(z)\right)^{1/2}<\infty,
\]
where $V$ denotes the Lebesgue measure on $\calU$. It is a closed subspace of $L^2(\calU)$ and hence a Hilbert space.
The orthogonal projection from $L^2(\calU)$ onto $A^2(\calU)$ can be expressed as an integral operator:
\[
Pf(z)=\int_{\calU} K(z,w)f(w)dV(w)
\]
with the Bergman kernel
\[
K(z,w)=\frac{n!}{4\pi^n}\left[ \frac{i}{2}(\overline{w}_n-z_n)-z^{\prime} \cdot \overline{w^{\prime}} \right]^{-n-1}.
\]
See, for instance, \cite[Theorem 5.1]{Gin64}. For simplicity, we write
\[
\bfrho(z,w)=\frac{i}{2}(\overline{w}_n-z_n)-z^{\prime} \cdot \overline{w^{\prime}}
\]
and let $\bfrho(z):=\bfrho(z,z)=\ImPt z_n-|z^{\prime}|^2$.
With this notation, the Bergman kernel of $\calU$ can be written as 
\[
K(z,w) = \frac {n!}{4\pi^n} \frac {1}{\bfrho(z,w)^{n+1}}, \quad z,w\in \calU.
\]

Given $f\in L^2(\calU)$, the Hankel operator with symbol $f$ is defined by
\[
H_f g=(I-P)(fg),\quad g\in A^2(\calU),
\]
where $I$ is the identity operator. In general, $H_f$ may be unbounded but it is densely defined by the fact that, for each $\alpha>n+1/2$, holomorphic functions $f$ on $\calU$
such that $f(z) = O(|z_n+i|^{-\alpha})$ form a dense subset of $A^2(\calU)$ (see Section 4).

Hankel operators have been widely studied by many authors in recent decades. The investigation of properties of Hankel operators was initiated by Axler \cite{Axl86} who characterized bounded and compact Hankel operators with analytic symbol on the Bergman spaces of the unit disk. That paper extended the classical results of Nehari \cite{Neh57}  and Hartman \cite{Har58} where the Hankel operators act on the Hardy space. Later on, Axler's result was generalized by Arazy  et al. \cite{AFP88,AFJP91} to weighted cases of the unit ball. For symbol functions bounded, Zhu \cite{Zhu87} first established the connection between size estimates of Hankel operators and the mean oscillation of the symbols in the Bergman metric. This idea was further investigated in a series papers \cite{BCZ87,BCZ88} and \cite{BBCZ90} in the context of bounded symmetric domains, and in \cite{Li92,Li94} in the context of bounded strongly pseudoconvex domains. In \cite{PZZ16}, Pau, Zhao and Zhu described the symbol functions such that Hankel operators are bounded or compact from $A^p_{\alpha}$ to $L^q_{\beta}$ on the unit ball, provided $1<p\leq q<\infty$. For recently popular Fock space, there has been a lot of researches on bounded and compact Hankel operators, see \cite{BY07,HW18,HL19,Sch04,SY13,WCZ15} for instance.
However, on the classical unbounded domains such as upper half-plane or upper half-space, over which there are few works on analogous results. See \cite{JP92} for a study of Hankel type operators on the upper half-plane.	

The main purpose of this paper is to characterize bounded and compact Hankel operators on $A^2(\calU)$. As in most above papers, our characterizations are based on the notions of mean oscillation and mean oscillation in the Bergman metric. 
Before stating our main results, we need to introduce some definitions and notations.

Let $b\calU:=\{z\in \mathbb{C}^n: \bfrho(z)=0\}$ denote the boundary of $\calU$. Then $\widehat{\calU}:=\calU\cup b\calU\cup\{\infty\}$ is the one-point compactification of $\calU$. Also, let $\partial \widehat{\calU}:=b\calU\cup\{\infty\}$. Thus, $z\to \partial\widehat{\calU}$ means $\bfrho(z)\to 0$ or $|z|\to\infty$.
We denote by $C_0(\calU)$ the space of complex-valued continuous functions $f$ on $\calU$ such that $f(z)\to 0$ as $z\to\partial \widehat{\calU}$.

Recall that the Berezin transform over $\calU$ is denoted by
\[
\widetilde{f}(z) = \int_{\calU} f(w)|k_z(w)|^2 dV(w),\quad z\in \calU,
\]
where
\[
k_z(w) = K(w,z)/\sqrt{K(z,z)},\quad w\in\calU.
\]
It is known from \cite{Liu2018} that the Berezin transform is bounded on $L^p(\calU)$ for $1<p\leq\infty$.
For $f\in L^2(\calU)$, the mean oscillation of $f$ at $z$ is denoted by
\[
MO(f)(z)=\left( \widetilde{|f|^2}(z)-|\widetilde{f}(z)|^2 \right)^{1/2}.
\]
It is easy to check that the function $MO(f)$ is positive and continuous. Let $BMO$ denote the space of functions $f\in L^2(\calU)$ such that 
\[
\|f\|_{BMO}= \sup\left\{MO(f)(z):z\in\calU\right\}
\]
is finite. We continue to say $f$ is in $VMO$ if $MO(f)$ is in $C_0(\calU)$.

The Bergman metric $\beta(\cdot,\cdot)$ is a complete Riemannian metric and is the ``integrated form'' of the infinitesimal metric 
\[
(g_{i,j}(z))= \frac {1}{n+1} \left(\frac {\partial^2 \log K(z,z)}{\partial z_i \partial \bar{z}_j}\right).
\]
We recall from \cite{LS20} that the Bergman metric on $\calU$ is given by 
\begin{equation*}\label{eqn:hyperdist}
 \beta(z,w)=\tanh^{-1}\sqrt{1-\frac{\bfrho(z)\bfrho(w)}{|\bfrho(z,w)|^2}}.
\end{equation*}
The Bergman metric ball at $z$ with radius $r>0$ is denoted by
 \[
 D(z,r) = \{w\in\calU: \beta(z,w)<r\}.
 \]
Let $|D(z,r)|=V(D(z,r))$. For fixed $r>0$ and $f\in L^2(\calU)$, we define the mean oscillation of $f$ at $z$ in the Bergman metric to be 
\[
MO_r(f)(z)= \left( \frac{1}{|D(z,r)|} \int_{D(z,r)} |f(w)-\widehat{f}_r(z)|^2 dV(w) \right)^{1/2},
\]
where 
\[
\widehat{f}_r(z)=\frac{1}{|D(z,r)|} \int_{D(z,r)} f(w) dV(w)
\]
is the integral mean of $f$ over $D(z,r)$. It is easy to check that the above two quantities are continuous functions on $\calU$. The space $BMO_r$ consists of those functions $f\in L^2(\calU)$ such that
\[
\|f\|_{BMO_r}= \sup\left\{MO_r(f)(z):z\in\calU\right\}
\]
is finite.
We continue to say $f$ is in $VMO_r$ if $MO_r(f)$ is in $C_0(\calU)$.

We require a few additional definitions. For $f$ continuous on $\calU$, we define 
\[
w(f)(z)=\sup\left\{|f(z)-f(w)|:w\in D(z,1)\right\}.
\]
Let $BO$ denote the space of continuous functions $f$ on $\calU$ such that $w(f)$ is bounded. We say $f$ is in $VO$ if $w(f)\in C_0(\calU)$. 

We denote by $BA$ the space of functions $f\in L^2(\calU)$ such that $\widetilde{|f|^2}$ is bounded,
or $VA$ if $f$ in addition satisfies $\widetilde{|f|^2}\in C_0(\calU)$.
Clearly, $BA\subset BMO$ and $VA\subset VMO$.

We now state our first two main results.

\begin{theorem}\label{thm:main1}
For $f\in L^2(\calU)$, the following conditions are equivalent:
\begin{enumerate}
  \item[(i)] $H_f$ and $H_{\overline{f}}$ are both bounded.
  \item[(ii)] $f\in BMO$.
  \item[(iii)] $f\in BMO_r$ for all $r>0$.
  \item[(iv)] $f\in BMO_r$ for some $r>0$.
  \item[(v)] $f\in BO+BA$.
\end{enumerate}
Moreover, the quantities $\max\{\| H_f\|,\|H_{\overline{f}}\|\}$, $\|f\|_{BMO}$ and $\|f\|_{BMO_r}$ are equivalent.
\end{theorem}

\begin{theorem}\label{thm:main2}
For $f\in L^2(\calU)$, the following conditions are equivalent:
\begin{enumerate}
  \item[(i)] $H_f$ and $H_{\overline{f}}$ are both compact.
  \item[(ii)] $f\in VMO$.
  \item[(iii)] $f\in VMO_r$ for all $r>0$.
  \item[(iv)] $f\in VMO_r$ for some $r>0$.
  \item[(v)] $f\in VO+VA$.
\end{enumerate}
\end{theorem}

We denote the Bloch space of $\calU$ by $\calB$, as the space of $f\in H(\calU)$ such that 
\[
\|f\|_{\calB}:=\sup\{|\widetilde{\nabla} f(z)|:z\in\calU\}
\]
is finite, where
\[
|\widetilde{\nabla}f(z)|=\left(2\bfrho(z)\left( 4\bfrho(z) \left|\frac{\partial f}{\partial z_n}(z)\right|^2+\sum_{j=1}^{n-1} \left|\frac{\partial f}{\partial z_j}(z) +2i\bar{z}_j \frac{\partial f}{\partial z_n}(z)\right|^2 \right) \right)^{1/2}.
\]
The operator $|\widetilde{\nabla}|$ is called invariant gradient of $\calU$, since it is M\"obius invariant, namely,  $|\widetilde{\nabla}(f\circ\psi)|=|\widetilde{\nabla}f\circ\psi|$ for all $f\in H(\calU)$ and all $\psi\in\mathrm{Aut}(\calU)$.
See Section 6 for details.
The little Bloch space of $\calU$ is denoted by $\calB_0$ consisting of functions $f\in\calB$ such that $|\widetilde{\nabla}f| \in C_0(\calU)$.

Our last main result is the following.

\begin{theorem}\label{thm:main3}
Let $H(\calU)$ be the space of holomorphic functions in $\calU$. For any $r>0$, we have 
\[
\calB=H(\calU)\cap BMO_r
\]
and 
\[
\calB_0=H(\calU)\cap VMO_r.
\]
Moreover, $\|f\|_{\calB}$ and $\|f\|_{BMO_r}$ are equivalent quantities.
\end{theorem}

For $f\in A^2(\calU)$, it is easy to see that $H_f=0$. In this case, as a consequence of Theorems \ref{thm:main1}, \ref{thm:main2} and \ref{thm:main3}, we obtain the following relationship between Bloch space and Hankel operators.

\begin{corollary}
For $f\in A^2(\calU)$, $H_{\overline{f}}$ is bounded if and only if $f\in\calB$, $H_{\overline{f}}$ is compact if and only if $f\in\calB_0$. Moreover, $\|H_{\overline{f}}\|$ and  $\|f\|_{\calB}$ are equivalent quantities.
\end{corollary}

We mention that, compared to the classical case of the unit disk or the unit ball, the results obtained in the above are not surprising and can be considered as the generalization of the classical case.  Nevertheless, the methods we used are different from the classical case in some respects, due to the unboundedness of the Siegel upper half-space, especially in the characterization of compact Hankel operators. In particular, the definition of Bloch space of the Siegel upper half-space has never been given in a direct way before.

This paper is organized as follows. Section 2 contains some auxiliary results. In Section 3, we study the relationships between $BMO$, $BMO_r$ and $VMO$, $VMO_r$, and as bridges, the spaces $BO$, $BA$ and $VO$, $VA$ are introduced and studied. These also show the equivalence of (ii)-(v) of both Theorems \ref{thm:main1} and \ref{thm:main2}, stated as Theorem \ref{thm:BMO,VMO}. Section 4 is devoted to characterizing bounded Hankel operators, where the equivalence of (i) and (ii) of Theorem \ref{thm:main1} is proved. In Section 5, we give a characterization of compact Hankel operators that proves the equivalence of (i) and (ii) of Theorem \ref{thm:main2}. In Section 6, we specialize to holomorphic functions and prove Theorem \ref{thm:main3}.

Throughout the paper, the letter $C$ will denote a positive constant that may vary at each occurrence but is independent of the essential variables.

\section{Preliminaries}

In this section we collect several basic lemmas.
For each $t>0$, we define the nonisotropic dilation $\delta_t$ by
\[
\delta_t(u)=(t u^{\prime},t^2 u_n), \quad u\in \calU.
\]
Also, to each fixed $z\in\calU$, we associate the following (holomorphic) affine
self-mapping of $\calU$:
\[
h_z(u) ~:=~ \left(u^{\prime} - z^{\prime}, u_n - \RePt z_n - 2i u^{\prime} \cdot \overline{z^{\prime}}  + i|z^{\prime}|^2 \right),
\quad u\in \calU.
\]
All these mappings are holomorphic automorphisms of $\calU$. See \cite[Chapter XII]{Ste93}. Hence the mappings 
\[
\sigma_z := \delta_{\bfrho(z)^{-1/2}} \circ h_z
\]
are holomorphic automorphisms of $\calU$. Simple calculations show that $\sigma_z(z)=\bfi:=(0^{\prime},i)$ and 
\begin{equation}\label{eqn:jacobian}
(J_R \sigma_z) (u) = \bfrho(z)^{-(n+1)},
\end{equation}
where $(J_R\sigma_z)(u)$ stands for the real Jacobian of $\sigma_z$ at $u$.
It is known that the Bergman metric is invariant under the holomorphic automorphisms of $\calU$; see for instance \cite[Proposition 1.4.15]{Kra01}. Hence, we have
\begin{equation}\label{eq:invarianceBeta}
\beta(\sigma_z(u),\sigma_z(v))=\beta(\sigma_z^{-1}(u),\sigma_z^{-1}(v))=\beta(u,v)
\end{equation}
for all $z,u,v\in\calU$.

\begin{lemma}\label{lem:rho(sigma_a)}
Suppose $z,u,v\in\calU$, we have 
\begin{align}
\bfrho(\sigma_z(u),\sigma_z(v))&= \bfrho(z)^{-1} \bfrho(u,v)\nonumber,\\
\bfrho(\sigma_z^{-1}(u),\sigma_z^{-1}(v))&= \bfrho(z) \bfrho(u,v)\label{eqn:-rho(sigma_a)}.
\end{align}
\end{lemma}

\begin{proof}
First note that
\[
\bfrho(\delta_t(u),\delta(v))=t^2 \bfrho(u,v).
\]
Also, an easy calculation shows that
\[
\bfrho(h_z(u),h_z(v))=\bfrho(u,v).
\]
Then a combination of the two above equalities gives
\begin{align*}
\bfrho(\sigma_z(u),\sigma_z(v)) &=\bfrho(\delta_{\bfrho(z)^{-1/2}}(h_z(u)),\delta_{\bfrho(z)^{-1/2}}(h_z(v)))\\
&=\bfrho(z)^{-1} \bfrho( h_z(u),h_z(v))\\
&= \bfrho(z)^{-1} \bfrho(u,v),
\end{align*}
which is the first equality of the lemma. 
The proof of the second one is exactly the same by observing that 
\[
\sigma_z^{-1}=h_z^{-1}\circ \delta_{\bfrho(z)^{1/2}}
\]
and
\[
\bfrho(h_z^{-1}(u),h_z^{-1}(v))=\bfrho(u,v).
\]
The proof is complete.
\end{proof}

The following three lemmas can be found in \cite{LS20}, they serve as basic tools in this paper.

\begin{lemma}\label{lem:elemtryeq1}
We have
\begin{equation}\label{eqn:elemtryeq1}
 2|\bfrho(z,w)|\geq \max\{\bfrho(z),\bfrho(w)\}
\end{equation}
for any $z, w\in \calU$.
\end{lemma}

\begin{lemma}\label{lem:volBergmanball}
For any $z\in\mathcal{U}$ and $r>0$ we have
\begin{equation}\label{eqn:volBergmanball}
|D(z,r)| ~=~ \frac {4\pi^n}{n!} \frac {\tanh^{2n} r } {(1-\tanh^2 r)^{n+1}}\, \bfrho(z)^{n+1} .
\end{equation}
\end{lemma}

\begin{lemma}\label{blem}
Given $r > 0$, the inequalities
\begin{equation}\label{eqn:eqvltquan}
\frac {1-\tanh (r)}{1+ \tanh (r)} ~\leq~ \frac{|\bfrho(z,u)|}{|\bfrho(z,v)|}
~\leq~ \frac {1+\tanh (r)}{1-\tanh (r)}
\end{equation}
hold for all $z,u,v\in \calU$ with $\beta(u,v)\leq r$.
\end{lemma}

\section{$BMO$, $VMO$, $BMO_r$ and $VMO_r$}

In this section we establish function-theoretic relationships between the spaces $BMO$, $BMO_r$ and $VMO$, $VMO_r$.  The main purpose of this section is to prove the following result that shows the equivalence of (ii)-(v) of both Theorems \ref{thm:main1} and \ref{thm:main2}.

\begin{theorem}\label{thm:BMO,VMO}
For any $r>0$, we have
$BMO=BO+BA=BMO_r$ and $VMO=VO+VA=VMO_r$.
Moreover, the quantities 
\[
\|f\|_{BMO},\quad \|f_1\|_{BO}+\|f_2\|_{BA},\quad \|f\|_{BMO_r},
\]
are equivalent for $f=f_1+f_2$ with $f_1\in BO$ and $f_2\in BA$.
\end{theorem}

To prove Theorem \ref{thm:BMO,VMO}, we first need the following equalities.
It is not hard to check that
\begin{align}
&MO_r(f)(z)=\left(\widehat{|f|^2_r}(z)-|\widehat{f_r}(z)|^2\right)^{1/2}\label{eq:MO_r1}\\
&=\left( \frac{1}{2|D(z,r)|^2} \int_{D(z.r)} \int_{D(z,r)} |f(u)-f(v)|^2 dV(u) dV(v)\right)^{1/2}\label{eq:MO_r2}.
\end{align}
It is also not hard to check that
\begin{align}
MO(f)(z)&=\left(\int_{\calU} |f(w)-\widetilde{f}(z)|^2 |k_z(w)|^2 dV(w)\right)^{1/2}\label{eq:MO1}\\
&=\left(\frac{1}{2} \int_{\calU} \int_{\calU} |f(u)-f(v)|^2 |k_z(u)|^2 |k_z(v)|^2 dV(u)dV(v)\right)^{1/2}\label{eq:MO2}.
\end{align}

\begin{lemma}\label{lem:BMOinBMO_r}
For any $r>0$, we have $BMO\subset BMO_r$ and $VMO\subset VMO_r$. Moreover, there is a positive constant such that $\|f\|_{BMO_r} \leq C \|f\|_{BMO}$.
\end{lemma}
\begin{proof}
By \eqref{eqn:volBergmanball} and \eqref{eqn:eqvltquan}, there is positive constant $C$ such that $|k_z(w)|^2\geq C V(D(z,r))^{-1}$ for $w\in D(z,r)$. Together with \eqref{eq:MO2} and \eqref{eq:MO_r2}, we have
\begin{align*}
MO(f)(z)&\geq \left(\frac{1}{2} \int_{D(z,r)} \int_{D(z,r)} |f(u)-f(v)|^2 |k_z(u)|^2 |k_z(v)|^2 dV(u)dV(v)\right)^{1/2}\\
&\geq \left(\frac{C^2}{2|D(z,r)|^2} \int_{D(z,r)} \int_{D(z,r)} |f(u)-f(v)|^2  dV(u) dV(v)\right)^{1/2}\\
&= C MO_r(f)(z).
\end{align*}
This proves the lemma.
\end{proof}

It remains to show that $BMO_r\subset BMO$ and $VMO_r\subset VMO$.
To this end, we need several auxiliary function spaces to link them up. 

For a continuous function $f$ on $\calU$, the oscillation $f$ at $z$ in the Bergman metric is denoted by 
\[
w_r(f)(z)=\sup\left\{|f(z)-f(w)|:w\in D(z,r)\right\}.
\]
Let $BO_r$ denote the space of continuous functions $f$ such that
\[
\|f\|_{BO_r}=\sup\left\{ w_r(f)(z):z\in\calU\right\}
\]
is finite. We continue to say $f$ is in $VO_r$ if $w_r(f)$ is in $C_0(\calU)$.

\begin{lemma}\label{lem:BO_r}
A continuous function $f$ belongs to $BO_r$ if and only if there exists a positive constant $C$ such that
\[
|f(z)-f(w)|\leq C (\beta(z,w) +1)
\]
for all $z,w\in\calU$.
As a consequence, the space $BO_r$ is independent of $r$.
\end{lemma}
\begin{proof}
The sufficiency is obviously true. It suffices to prove the necessity. Suppose $f\in BO_r$. If $\beta(z,w)\leq r$, then the desired inequality is obvious. So fix any $\beta(z,w)>r$, and let $\gamma(t)$, $0\leq t\leq 1$, be the geodesic in the Bergman metric from $z$ to $w$. Write $N=[r^{-1} \beta(z,w)]+1$, where $[x]$ denotes the greatest integer less than or equal to $x$. Let $t_i = i/N$, $0\leq i\leq N$. Then
\[
\beta(\gamma(t_{i}),\gamma(t_{i-1}))=\frac{\beta(z,w)}{N} \leq r.
\]
Since $f\in BO_r$, there is a constant $M>0$ such that 
\[
|f(u)-f(v)|\leq M
\]
whenever $\beta(u,v)\leq r$. Thus, 
\begin{align*}
|f(z)-f(w)| &\leq \sum_{i=1}^N |f(\gamma(t_{i}))-f(\gamma(t_{i-1}))|\\
&\leq N M \leq M (1+r^{-1} \beta(z,w))\\
&\leq M (1+r^{-1}) (\beta(z,w)+1).
\end{align*}
This proves the desired inequality.
\end{proof}

\begin{lemma}\label{lem:VO_r}
The space $VO_r$ is independent of $r$.
\end{lemma}
\begin{proof}
Follow the proof of Lemma \ref{lem:BO_r}.
\end{proof}

Since the spaces $BO_r$ and $VO_r$ are both independent of $r$, we will simply regard $BO_r$ and $VO_r$ as $BO$ and $VO$ respectively. Furthermore, we write
\[
\|f\|_{BO}=\sup\left\{ |f(z)-f(w)|:\beta(z,w)\leq 1\right\}.
\]
The quantities $\|f\|_{BO}$ and $\|f\|_{BO_r}$ are equivalent for any $r>0$.

Let $BA_r$ denote the space of locally Lebesgue integral functions $f$ such that
\[
\|f\|_{BA_r}=\sup\left\{ \left[\widehat{|f|^2_r}(z)\right]^{1/2}:z\in\calU \right\}
\]
is finite. We continue to say $f$ is in $VA_r$ if $\left[\widehat{|f|^2_r}(\cdot)\right]^{1/2}$ is in $C_0(\calU)$.
Recall from \cite{LS20} that $f\in BA_r$ if and only if $\widetilde{|f|^2}$ is bounded,
$f\in VA_r$ if and only if $\widetilde{|f|^2}$ is in $C_0(\calU)$. 
Moreover, the quantities $\|f\|_{BA_r}$ and
\[
\|f\|_{BA}=\sup \left\{ \left[ \widetilde{|f|^2}(z) \right]^{1/2} : z\in\calU \right\}
\]
are equivalent for any $r>0$.
It follows that the spaces $BA_r$ and $VA_r$ are both independent of $r$. So we will simply write $BA$ for $BA_r$ and $VA$ for $VA_r$.

\begin{lemma}\label{lem:normBOA}
Suppose $r>0$ and $f\in BMO_{2r}$. Then $\widehat{f_r}\in BO$ and $f-\widehat{f_r}\in BA$ with
\[
\|\widehat{f_r}\|_{BO_r}\leq C\|f\|_{BMO_{2r}}, \quad \|f-\widehat{f_r}\|_{BA_r}\leq C\|f\|_{BMO_{2r}}
\]
for some positive constant $C$. 
Furthermore, for any $r>0$, we have $BMO_{r}\subset BO+BA$ and 
\[
\|f_1\|_{BO}+\|f_2\|_{BA} \leq C \|f\|_{BMO_{r}}
\]
for $f=f_1+f_2$ with $f_1\in BO$ and $f_2\in BA$.
\end{lemma}
\begin{proof}
Given $z,w\in\calU$ with $\beta(z,w)<r$, the Cauchy-Schwarz inequality yields
\[
|\widehat{f_r}(z)-\widehat{f_r}(w)|^2\leq \frac{1}{|D(z,r)| |D(w,r)|}\int_{D(z.r)} \int_{D(w,r)} |f(u)-f(v)|^2 dV(u) dV(v),
\]
which, by \eqref{eqn:volBergmanball} and \eqref{eqn:eqvltquan}, is no greater than
\[
\frac{C}{|D(z,r)|^2} \int_{D(z.2r)} \int_{D(z,2r)} |f(u)-f(v)|^2 dV(u) dV(v)
\]
for a suitable positive constant $C$. In view of \eqref{eq:MO_r2}, the above quantity is equal to $2CMO_{2r}^2(f)(z)$. It follows that 
\begin{equation}\label{eq:normBOA1}
w_r(\widehat{f_r})(z)\leq \sqrt{2C} MO_{2r}(f)(z)
\end{equation}
and hence
\[
\|\widehat{f_r}\|_{BO_r}\leq \sqrt{2C}\|f\|_{BMO_{2r}}.
\]
Next, we consider $g=f-\widehat{f_r}$. By the triangle inequality we have
\begin{align*}
\left[\widehat{|g|^2_r}(z)\right]^{1/2}
&= \left( \frac{1}{|D(z,r)|} \int_{D(z,r)} |f(w)-\widehat{f}_r(w)|^2 dV(w) \right)^{1/2}\\
&\leq \left( \frac{1}{|D(z,r)|} \int_{D(z,r)} |f(w)-\widehat{f}_r(z)|^2 dV(w) \right)^{1/2}\\
&\quad + \left( \frac{1}{|D(z,r)|} \int_{D(z,r)} |\widehat{f}_r(z)-\widehat{f}_r(w)|^2 dV(w) \right)^{1/2}\\
&\leq MO_r(f)(z)+w_r(\widehat{f_r})(z).
\end{align*}
By \eqref{eq:MO_r2} and \eqref{eqn:volBergmanball}, there is another positive constant $C_1$ such that 
\[
MO_r(f)(z)\leq C_1 MO_{2r}(f)(z).
\]
This together with \eqref{eq:normBOA1} yields that
\begin{equation}\label{eq:normBOA2}
\left[\widehat{|g|^2_r}(z)\right]^{1/2} \leq (\sqrt{2C}+C_1) MO_{2r}(f)(z)
\end{equation}
and hence
\[
\|f-\widehat{f_r}\|_{BA_r}\leq (\sqrt{2C}+C_1) \|f\|_{BMO_{2r}}.
\]
This completes the proof of the lemma.
\end{proof}

By \eqref{eq:normBOA1} and \eqref{eq:normBOA2}, we also obtain the following result.
\begin{lemma}\label{lem:VOA}
Suppose $r>0$ and $f\in VMO_{2r}$. Then $\widehat{f_r}\in VO$ and $f-\widehat{f_r}\in VA$. Furthermore, for any $r>0$, we have $VMO_{r}\subset VO+VA$ for $f=f_1+f_2$ with $f_1\in VO$ and $f_2\in VA$.
\end{lemma}



Now the aim of the section will be achieved if we can show that $BO+BA\subset BMO$ and $VO+VA\subset VMO$. First, note that 
\[
MO(f)(z)^2=\widetilde{|f|^2}(z)-|\widetilde{f}(z)|^2 \leq \widetilde{|f|^2}(z).
\]
This immediately yields the following result.

\begin{lemma}
We have $BA\subset BMO$ and $VA\subset VMO$. Moreover, $\|f\|_{BMO}\leq \|f\|_{BA}$.
\end{lemma}

Next, we shall prove that $BO\subset BMO$ and $VO\subset VMO$. For this, we need characterizations for $BO$ and $VO$. 

\begin{lemma}\label{lem:thm:chaBO}
If $f\in BMO$, then $\widetilde{f}-\widehat{f_r}$ is bounded continuous and $\widetilde{f}$ is in $BO$. 
If $f\in VMO$, then $\widetilde{f}-\widehat{f_r}$ is in $C_0(\calU)$ and $\widetilde{f}$ is in $VO$. 
\end{lemma}
\begin{proof}
Suppose $f\in BMO$ (or $VMO$). Then by Lemma \ref{lem:BMOinBMO_r}, $f\in BMO_{2r}$ (or $VMO_{2r}$). This along with Lemma \ref{lem:normBOA} (or Lemma \ref{lem:VOA}) implies that $\widehat{f_r}\in BO$ (or $VO$).

Since
\[
|\widehat{f_r}(z)-\widetilde{f}(z)|\leq \frac{1}{|D(z,r)|} \int_{D(z,r)} |f(w)-\widetilde{f}(z)| dV(w),
\]
by Cauchy-Schwarz inequality and a combination of  \eqref{eqn:volBergmanball} and \eqref{eqn:eqvltquan}  there exists a positive constant $C$ such that
\begin{align*}
|\widehat{f_r}(z)-\widetilde{f}(z)|^2 &\leq \frac{1}{|D(z,r)|} \int_{D(z,r)} |f(w)-\widetilde{f}(z)|^2 dV(w)\\
& \leq C \int_{D(z,r)} |f(w)-\widetilde{f}(z)|^2 |k_z(w)|^2 dV(w)\\
&\leq C MO(f)(z)^2
\end{align*}
in view of \eqref{eq:MO1}. It follows that $\widetilde{f}-\widehat{f_r}$ is bounded continuous (or in $C_0(\calU)$). Since $\widehat{f_r}\in BO$ (or $VO$), $\widetilde{f}\in BO$ (or $VO$). The proof of the lemma is complete.
\end{proof}

\begin{lemma}\label{lem:simkeylem}
Let $t>-1$, $s-t>n+1$ and $\alpha\geq 0$. There is a positive constant $C$ such that 
\[
\int_{\calU} \frac{\beta(z,w)^{\alpha} \bfrho(w)^t}{|\bfrho(z,w)|^s} dV(w)\leq C \bfrho(z)^{n+1+t-s}
\]
for all $z\in\calU$.
\end{lemma}
\begin{proof}
Since $\log x<x^{\epsilon}$ holds for any $x>0$ and any $\epsilon>0$, it follows that
\[
\beta(z,w)\leq \frac{1}{2}\log\frac{4|\bfrho(z,w)|^2}{\bfrho(z)\bfrho(w)}\leq 2^{2\epsilon-1}
\frac{|\bfrho(z,w)|^{2\epsilon}}{\bfrho(z)^{\epsilon}\bfrho(w)^{\epsilon}}.
\]
Thus, the desired integral
\[
\int_{\calU} \frac{\beta(z,w)^{\alpha} \bfrho(w)^t}{|\bfrho(z,w)|^s} dV(w)
\leq 2^{(2\epsilon-1)\alpha} \bfrho(z)^{-\epsilon\alpha} 
\int_{\calU} \frac{\bfrho(w)^{t-\epsilon\alpha}}{|\bfrho(z,w)|^{s-2\epsilon\alpha}} dV(w).
\]
Let $\epsilon$ be small enough such that $t-\epsilon\alpha>-1$ and $s-t-\epsilon\alpha>n+1$. Then by \cite[Lemma 5]{LLHZ18}  there is a positive constant $C$ such that
\[
\int_{\calU} \frac{\bfrho(w)^{t-\epsilon\alpha}}{|\bfrho(z,w)|^{s-2\epsilon\alpha}} dV(w)
= C \bfrho(z)^{n+1+t+\epsilon\alpha-s}.
\]
This completes the proof of the lemma. 
\end{proof}

\begin{lemma}\label{thm:chaBO}
For a continuous function $f$ on $\calU$, the following conditions are equivalent:
\begin{enumerate}
  \item [(i)] $f\in BO$.
  \item [(ii)] There is a positive constant $C$ such that 
  \[
  |f(z)-f(w)|\leq C (\beta(z,w) +1)
  \]
  for all $z,w\in\calU$.
  \item [(iii)] $\|(f(z)-f)k_z\|$ is a bounded continuous  function on $\calU$.
  \item [(iv)] $f\in BMO$ and $f-\widetilde{f}$ is a bounded continuous function on $\calU$.
\end{enumerate}
Moreover, we have $\|f\|_{BMO}\leq C\|f\|_{BO}$.
\end{lemma}
\begin{proof}
(i)$\Rightarrow$(ii) follows from Lemma \ref{lem:BO_r}. For (ii)$\Rightarrow$(iii), by making a change of variables with $\sigma_z^{-1}$, then \eqref{eqn:jacobian} and \eqref{eqn:-rho(sigma_a)}, we have
\begin{align*}
\int_{\calU}  |f(z)-f(w)|^2 |k_z(w)|^2 dV(w)
= \int_{\calU} |f(\sigma_z^{-1}(\bfi))-f(\sigma_z^{-1}(u))|^2 |k_{\bfi}(u)|^2 dV(u),
\end{align*}
which, by assumption, is less than or equal to 
\[
C^2  \int_{\calU} [1+\beta(\sigma_z^{-1}(\bfi),\sigma_z^{-1}(u))]^2 |k_{\bfi}(u)|^2 dV(u).
\]
In view of \eqref{eq:invarianceBeta}, the proceeding expression is equal to 
\[
C^2  \int_{\calU} [1+\beta(\bfi,u)]^2 |k_{\bfi}(u)|^2 dV(u).
\]
It follows from Lemma \ref{lem:simkeylem} that the above integral is bounded and hence (iii) follows.

(iii)$\Rightarrow$(iv) follows by observing that 
\begin{equation}\label{eq:thm:chaBO}
\|(f(z)-f)k_z\|^2=MO(f)(z)^2 + |f(z)-\widetilde{f}(z)|^2.
\end{equation}

Finally, assume that $f\in BMO$ and $f-\widetilde{f}$ is bounded continuous. By the first assertion of Lemma \ref{lem:thm:chaBO}, we see that $\widetilde{f}\in BO$. It follows that $f\in BO$, so (iv)$\Rightarrow$(i) follows.
\end{proof}

\begin{lemma}\label{thm:chaVO}
For a continuous function $f$ on $\calU$, the following conditions are equivalent:
\begin{enumerate}
  \item [(i)] $f\in VO$.
  \item [(ii)] $\|(f(z)-f)k_z\|\in C_0(\calU)$.
  \item [(iii)] $f\in VMO$ and $f-\widetilde{f}\in C_0(\calU)$.
\end{enumerate}
\end{lemma}
\begin{proof}
(i)$\Rightarrow$(ii). On one hand, reviewing the proof of Theorem \ref{thm:chaBO}, we see that
\[
\|(f(z)-f)k_z\|=\|(f(z)-f\circ\sigma_z^{-1})k_{\bfi}\|,
\]
and for each $z\in\calU$, the function $(f(z)-f\circ\sigma_z^{-1}(\cdot))k_{\bfi}(\cdot)$ is uniformly bounded by the square integrable function
$(1+\beta(\bfi,\cdot))k_{\bfi}(\cdot)$. On the other hand, it follows from Lemma \ref{lem:VO_r} that
\[
f(z)-f\circ\sigma_z^{-1}(w)\to 0 \quad \text{as} \quad z\to \partial\widehat{\calU}
\]
for each $w\in\calU$ (since $\beta(z,\sigma_z^{-1}(w))=\beta(\bfi,w)$). This implies that
\[
(f(z)-f\circ\sigma_z^{-1}(w))k_{\bfi}(w) \to 0 \quad \text{as} \quad z\to \partial\widehat{\calU}
\]
for each $w\in\calU$. Therefore, we conclude by the Lebesgue dominated convergence theorem that
$\|(f(z)-f\circ\sigma_z^{-1})k_{\bfi}\|$ is in $C_0(\calU)$ and so is $\|(f(z)-f)k_z\|$.

That implication (ii)$\Rightarrow$(iii) follows by \eqref{eq:thm:chaBO}.

(iii)$\Rightarrow$(i). Assume (iii). It follows from the second assertion of Lemma \ref{lem:thm:chaBO} that $\widetilde{f}\in VO$ and hence $f\in VO$. The proof is complete.
\end{proof}

It follows from the above two lemmas that $BO\subset BMO$ and $VO\subset VMO$. Therefore, the proof Theorem \ref{thm:BMO,VMO} is complete.


\section{Bounded Hankel operators}\label{sec:boundedHankel}

The main purpose of this section is to prove the following theorem that proves the equivalence of (i) and (ii) of Theorem \ref{thm:main1}.

\begin{theorem}\label{thm:boundHankel}
Let $f\in L^2(\calU)$. Then $f$ is in $BMO$ if and only if both $H_f$ and $H_{\overline{f}}$ are bounded on $A^2(\calU)$. Moreover, there is a positive constant $C$ such that 
\[
C^{-1} \|f\|_{BMO} \leq \|H_f\| + \|H_{\overline{f}}\| \leq C \|f\|_{BMO}
\]
for all $f\in BMO$.
\end{theorem}

First of all, we need to check that Hankel operators $H_f$ with symbol $f\in L^2(\calU)$ are well defined on $A^2(\calU)$. Given $f\in L^2(\calU)$, apparently, we only need to show that the multiplication operator $M_f$ given by $M_fg=fg$ is dense defined from $A^2(\calU)$ to $L^2(\calU)$. 
Recall  from \cite[Theorem 4.1]{LS20} that the space denoted by $\calS_{\alpha}$ consisting of holomorphic functions $g$ in $\calU$ such that 
\[
\sup_{z\in\calU} |z_n+i|^{\alpha} |g(z)| <\infty
\]
is dense in $A^2(\calU)$ whenever $\alpha>n+1/2$. 
Note that $|z_n+i|\geq 1$ for all $z\in\calU$. It follows that the functions in $\calS_{\alpha}$ with $\alpha>n+1/2$ are bounded on $\calU$. Therefore, $M_f$ is dense defined from $A^2(\calU)$ to $L^2(\calU)$. So $H_f$ is well defined on $A^2(\calU)$ as desired.

Now we first prove the necessity of Theorem \ref{thm:boundHankel}.

\begin{proposition}
There is a positive constant $C$ such that $\|H_f\|\leq C\|f\|_{BO}$ for all $f\in BO$.
\end{proposition}
\begin{proof}
Suppose $f\in BO$. It follows from Lemma \ref{thm:chaBO} that there exists a positive constant $C$ such that 
\[
|f(z)-f(w)|\leq C \|f\|_{BO} (\beta(z,w)+1)
\]
for all $z,w\in\calU$. Recalling the integral form of Hankel operators that
\[
H_f g(z) = \int_{\calU} f(z)-f(w) g(w) K(z,w) dV(w),
\]
we have
\begin{align*}
|H_f g(z)| &\leq \int_{\calU} |f(z)-f(w)| |g(w)| |K(z,w)| dV(w)\\
&\leq C \frac{n!}{4\pi^n} \|f\|_{BO} \int_{\calU} \frac{\beta(z,w)+1}{|\bfrho(z,w)|^{n+1}} |g(w)| dV(w)
\end{align*}
for $g\in A^2(\calU)$. It follows from Theorems 1 and 6 of \cite{LLHZ18} that the operator
\[
g \mapsto \int_{\calU} \frac{\beta(\cdot,w)+1}{|\bfrho(\cdot,w)|^{n+1}} |g(w)| dV(w)
\]
is bounded on $L^2(\calU)$. Thus, the desired result follows.
\end{proof}


\begin{proposition}\label{prop:HankelBA}
There is a positive constant $C>0$ such that $\|H_f\|\leq C\|f\|_{BA}$ for all $f\in BA$.
\end{proposition}
\begin{proof}
Suppose $f\in BA$. It follows from \cite[Theorem 1.1]{LS20} that the operator $i_f: A^2(\calU)\to L^2(\calU,d\mu_f)$ is bounded and $\|i_f\|$ is comparable to $\|f\|_{BA}$. Here $d\mu_f=|f|^2 dV$.
Let $g\in\calS_{\alpha}$ with $\alpha>n+1/2$, which is dense in $A^2(\calU)$. Then
\begin{align*}
\|H_f g\|^2 &=\|(I-P)(fg)\|^2 \leq \|fg\|^2=\int_{\calU} |g(z)|^2 d\mu_f(z)\\
&\leq C \|f\|_{BA} \int_{\calU} |g(z)|^2 dV(z)
\end{align*}
for some positive constant $C$. This proves the lemma.
\end{proof}

Combining the above two propositions with the decomposition
\[
BMO=BO+BA
\]
given in Theorem \ref{thm:BMO,VMO}, we see that if $f\in BMO$, then $H_f$ and $H_{\overline{f}}$ are both bounded on $A^2(\calU)$ with 
\[
\|H_f\|\leq C \|f\|_{BMO}, \quad\quad \|H_{\overline{f}}\|\leq C \|f\|_{BMO}
\]
for a suitable positive constant $C$ independent of $f$. Therefore, the necessity of Theorem \ref{thm:boundHankel} follows.

Next, we turn to prove the sufficiency of Theorem \ref{thm:boundHankel}. We should point out that the idea of this part is inspired by that of Pau, Zhao and Zhu \cite{PZZ16}.

\begin{lemma}\label{lem:lu_suf_thm_boundHankel}
For every $z\in\calU$ and $f\in L^2(\calU)$, the function $f_z$ given by
\[
f_z(w)=\frac{P(\overline{f} k_z)(w)}{k_z(w)}, \quad w\in\calU,
\]
is holomorphic and $P(\overline{f_z} k_z)= \widetilde{f}(z) k_z$.
\end{lemma}
\begin{proof}
Fix $z\in\calU$ and $f\in L^2(\calU)$. By \eqref{eqn:elemtryeq1} we have 
\[
|K(z,w)|\leq 2^{n+1} K(z,z).
\]
Thus, we have
\[
\int_{\calU} |f(w) k_z(w)|^2 dV(w)\leq 4^{n+1} K(z,z) \int_{\calU} |f(w)|^2 dV(w),
\]
which implies that $\overline{f}k_z\in L^2(\calU)$. Also, since $k_z$ never vanishes on $\calU$, the function $f_z$ is well defined and holomorphic. 

To prove the second assertion, we first need to verify that $P(\overline{f_z} k_z)$ is well defined.
In fact, since $f_z k_z = P(\overline{f} k_z)$ and $\overline{f}k_z\in L^2(\calU)$, by the boundedness of $P$ we have $f_z k_z\in A^2(\calU)$. As a consequence, $P(\overline{f_z} k_z)$ is well defined. 

Claim: If $K(\cdot,w) h(\cdot) \in A^2(\calU)$ for every $w\in\calU$, then $P(\overline{h} k_z)=\overline{h(z)} k_z$. It follows from the reproducing formula of functions in Bergman spaces on $\calU$ (see \cite{DK93} for instance) that 
\begin{align*}
P(\overline{h} k_z)(w)&=\int_{\calU} K(w,u)\overline{h(u)} k_z(u) dV(u)\\
&=\frac{1}{\sqrt{K(z,z)}} \overline{\int_{\calU} K(u,w) h(u) K(z,u) dV(u)}\\
&= \frac{\overline{K(z,w) h(z)}}{\sqrt{K(z,z)}}=\overline{h(z)} k_z(w),
\end{align*}
as claimed.

For any fixed $w\in\calU$, notice from \eqref{eqn:eqvltquan} that the ratio $|K(u,w)/K(u,z)|$ is bounded for all $u\in\calU$. This together with the fact $f_z k_z\in A^2(\calU)$ implies that $f_z(\cdot) K(\cdot,w)\in A^2(\calU)$ for every $w\in\calU$. Thus, by the claim we obtain $P(\overline{f_z} k_z)=\overline{f_z(z)} k_z$. 
Also, note that 
\begin{align*}
\overline{f_z(z)}&=\frac{\overline{P(\overline{f}k_z)(z)}}{\sqrt{K(z,z)}}=\frac{1}{\sqrt{K(z,z)}}\overline{\int_{\calU} K(z,u) \overline{f(u)} k_z(u) dV(u)}\\
&=\frac{1}{\sqrt{K(z,z)}}\int_{\calU}K(u,z) f(u) \overline{k_z(u)} dV(u)\\
&=\int_{\calU} f(u) |k_z(u)|^2 dV(u) = \widetilde{f}(z).
\end{align*}
Hence, 
\[
P(\overline{f_z} k_z)(w)=\widetilde{f}(z) k_z(w)
\]
as desired.
\end{proof}

\begin{lemma}\label{lem:suf_thm_boundHankel}
There is a positive constant $C$ such that 
\[
MO(f)(z)\leq C \left( \|H_fk_z\| + \|H_{\overline{f}}k_z\|\right)
\]
for all $f\in L^2(\calU)$ and all $z\in\calU$.
\end{lemma}
\begin{proof}
By Lemma \ref{lem:lu_suf_thm_boundHankel}, we know that $P(\overline{f_z} k_z)= \widetilde{f}(z) k_z$.
Recalling \eqref{eq:MO1}, we have
\[
MO(f)(z)=\|fk_z-\widetilde{f}(z)k_z\|=\|fk_z-P(\overline{f_z} k_z)\|.
\]
By the triangle inequality and the definition of Hankel operators, we obtain
\begin{align*}
MO(f)(z) &\leq \|fk_z-P(fk_z)\|+\|P(fk_z)-P(\overline{f_z} k_z)\| \\
&= \|H_f k_z\|+\|P(fk_z)-P(\overline{f_z} k_z)\|.
\end{align*}
It suffices to estimate the second term. Observe that
\begin{align*}
\|P(fk_z)-P(\overline{f_z} k_z)\| &\leq \|P\| \|fk_z-\overline{f_z} k_z\| 
= \|P\| \|\overline{f}k_z -f_z k_z\| \\
&= \|P\| \|\overline{f}k_z-P(\overline{f}k_z)\|=\|P\| \|H_{\overline{f}}k_z\|.
\end{align*}
This proves the result with constant $C=1+\|P\|$.
\end{proof}

Since each $k_z$ is a unit vector in $A^2(\calU)$, we see that if both $H_f$ and $H_{\overline{f}}$ are bounded on $A^2(\calU)$, then $f\in BMO$ with 
\[
\|f\|_{BMO} \leq C \left(\|H_f\|+\|H_{\overline{f}}\|\right),
\]
where $C$ is a positive constant independent of $f$. Thus, the sufficiency of Theorem \ref{thm:boundHankel} is proved.

\section{Compact Hankel operators}

The main result of this section is the following characterization of compact Hankel operators that proves the equivalence of (i) and (ii) of Theorem \ref{thm:main2}.

\begin{theorem}\label{thm:compactHankel}
Let $f\in L^2(\calU)$. Then $f$ is in $VMO$ if and only if both $H_f$ and $H_{\overline{f}}$ are compact on $A^2(\calU)$. 
\end{theorem}

Recall from \cite[Lemma 2.12]{LS20} that $k_z\to 0$ weakly in $A^2(\calU)$ as $z\to\partial\widehat{\calU}$, then the sufficiency of Theorem \ref{thm:compactHankel} is clear by using Lemma \ref{lem:suf_thm_boundHankel}. Thus, it suffices to prove the necessity. 

\begin{lemma}\label{lem:VAHankel}
If $f\in VA$, then $H_f$ is compact.
\end{lemma}
\begin{proof}
Suppose $f\in VA$ and $\{g_k\}$ is a sequence in $A^2(\calU)$ that converges to 0 weakly. 
Following the proof of Proposition \ref{prop:HankelBA}, we have
\[
\|H_f g_k\|^2 \leq \int_{\calU} |g_k(z)|^2 d\mu_f(z),
\]
where $d\mu_f=|f|^2 dV$. Since $f\in VA$, it follows from \cite[Theorem 1.2]{LS20} that the operator $i_f: A^2(\calU)\to L^2(\calU,d\mu_f)$ is compact.
The desired result then follows.
\end{proof}

\begin{lemma}\label{lem:deriveVO}
Let $0<\gamma<1$ and $\alpha\geq 0$. If $f\in VO$, then 
\[
\bfrho(z)^{\gamma} \int_{\calU} \frac{|f(z)-f(w)|^{\alpha} \bfrho(w)^{-\gamma}}{|\bfrho(z,w)|^{n+1}} dV(w)\to 0
\] 
as $z\to\partial\widehat{\calU}$.
\end{lemma}
\begin{proof}
By making a change of variables, \eqref{eqn:jacobian} and \eqref{eqn:-rho(sigma_a)}, we have
\begin{align*}
&\bfrho(z)^{\gamma} \int_{\calU} \frac{|f(z)-f(w)|^{\alpha} \bfrho(w)^{-\gamma}}{|\bfrho(z,w)|^{n+1}} dV(w) \\
&=\bfrho(z)^{\gamma} \int_{\calU} \frac{|f(z)-f(\sigma_z^{-1}(u))|^{\alpha} \bfrho(\sigma_z^{-1}(u))^{-\gamma}}{|\bfrho(\sigma_z^{-1}(\bfi),\sigma_z^{-1}(u))|^{n+1}} \bfrho(z)^{n+1} dV(u)\\
&=\int_{\calU} \frac{|f(z)-f(\sigma_z^{-1}(u))|^{\alpha}\bfrho(u)^{-\gamma}}{|\bfrho(\bfi,u)|^{n+1}} dV(u):=\int_{\calU} g_z(u) dV(u).
\end{align*}
On one hand, since $f\in VO$, from the proof of (i) $\Rightarrow$ (ii) of Theorem \ref{thm:chaVO} we see that $g_z(u)\to 0$ as $z\to\partial\widehat{\calU}$ for each fixed $u\in\calU$. On the other hand, it follows from the proof of Theorem \ref{thm:chaBO} that $g_z(u)$ is bounded by an integrable function 
\[
\frac{(1+\beta(\bfi,u))^{\alpha}\bfrho(u)^{-\gamma}}{|\bfrho(\bfi,u)|^{n+1}},
\]
of which the integrability is guaranteed by Lemma \ref{lem:simkeylem}. Therefore, it follows from the Lebesgue dominated convergence theorem that 
\[
\int_{\calU} g_z(u) dV(u)\to 0 \quad \text{as} \quad z\to\partial\widehat{\calU},
\]
completing the proof of the lemma.
\end{proof}

\begin{corollary}\label{lem:anopropVO}
Let $0<\gamma<1$ and $\alpha\geq 0$. If $f\in VO$, we have
\[
\lim_{R\to\infty} \sup_{z\in D(\bfi,R)^c} 
\bfrho(z)^{\gamma} \int_{\calU} \frac{|f(z)-f(w)|^{\alpha} \bfrho(w)^{-\gamma}}{|\bfrho(z,w)|^{n+1}} dV(w)= 0.
\]
\end{corollary}
\begin{proof}
By Lemma \ref{lem:deriveVO} and the definition of $\partial\widehat{\calU}$, it suffices to prove that $\beta(z,\bfi)\to\infty$ if and only if $\bfrho(z)\to 0$ or $|z|\to\infty$.
Since
\[
\beta(z,\bfi)=\tanh^{-1}\sqrt{1-\bfrho(z)/|\bfrho(z,\bfi)|^2},
\]
it is obvious that $\beta(z,\bfi)\to\infty$ if and only if $\bfrho(z)/|\bfrho(z,\bfi)|^2\to 0$.
Thus, it reduces to prove that $\bfrho(z)\to 0$ or $|z|\to\infty$ if and only if $\bfrho(z)/|\bfrho(z,\bfi)|^2\to 0$.
Note that, for $z\in\calU$, we have $1\leq |z_n+i|=2|\bfrho(z,\bfi)| \leq |z|+1$ and $\bfrho(z)\leq |z|\leq 2|\bfrho(z,\bfi)|$.
Therefore, we deduce that
\[
4\bfrho(z)/(|z|+1)^2 \leq \bfrho(z)/|\bfrho(z,\bfi)|^2 \leq 4\min\{\bfrho(z),|z|^{-1} \}
\]
for $z\in\calU$. Then  the lemma follows. 
\end{proof}

\begin{lemma}\label{lem:VOHankel}
If $f\in VO$, then $H_f$ is compact.
\end{lemma}
\begin{proof}
Suppose $f\in VO$. Let $\{g_k\}$ be a sequence in $A^2(\calU)$ that converges to 0 weakly. Then by \cite[Lemma 2.11]{LS20}, $\{g_k\}$ is bounded in $A^2(\calU)$ and converges to 0 uniformly on each compact subset of $\calU$. 

Let $0<\gamma<1/2$. For any  given $\varepsilon>0$, it follows from Corollary \ref{lem:anopropVO} that there exsists a sufficiently large $R$ such that
\begin{equation}\label{eq:lem:VOHankel}
\sup_{w\in D(\bfi,R)^c}\bfrho(w)^{2\gamma} \int_{\calU} \frac{|f(z)-f(w)|^2 \bfrho(z)^{-2\gamma}}{|\bfrho(z,w)|^{n+1}} dV(z)\leq \varepsilon^2.
\end{equation}
Recalling the integral form of Hankel operators, we have
\[
\|H_f g_k\| \leq \left[ \int_{\calU} \left( \int_{\calU} |f(z)-f(w)| |K(z,w)| |g_k(w)| dV(w) \right)^2 dV(z) \right]^{1/2}.
\]
We split the above integral into two parts:
\[
I_1=\left[ \int_{\calU} \left( \int_{D(\bfi,R)^c} |f(z)-f(w)| |K(z,w)| |g_k(w)| dV(w) \right)^2 dV(z) \right]^{1/2}
\]
and 
\[
I_2=\left[ \int_{\calU} \left( \int_{D(\bfi,R)} |f(z)-f(w)| |K(z,w)| |g_k(w)| dV(w) \right)^2 dV(z) \right]^{1/2}.
\]
We first deal with $I_1$.
By Cauchy-Schwarz inequality and Lemma \ref{lem:simkeylem}, the inner integral of $I_1$ is dominated by
\begin{align*}
&\left( \int_{D(\bfi,R)^c} \frac{|f(z)-f(w)| |g_k(w)| \bfrho(w)^{\gamma}}{|\bfrho(z,w)|^{(n+1)/2}} 
\frac{\bfrho(w)^{-\gamma}}{|\bfrho(z,w)|^{(n+1)/2}}dV(w) \right)^2\\
&\leq \int_{D(\bfi,R)^c} \frac{|f(z)-f(w)|^2 |g_k(w)|^2 \bfrho(w)^{2\gamma}}{|\bfrho(z,w)|^{n+1}} dV(w)\\
&\quad \times \int_{D(\bfi,R)^c} \frac{\bfrho(w)^{-2\gamma}}{|\bfrho(z,w)|^{n+1}} dV(w)\\
&\leq C\bfrho(z)^{-2\gamma} \int_{D(\bfi,R)^c} \frac{|f(z)-f(w)|^2 |g_k(w)|^2 \bfrho(w)^{2\gamma}}{|\bfrho(z,w)|^{n+1}} dV(w)
\end{align*}
for some positive constant $C$.
This together with \eqref{eq:lem:VOHankel} gives
\begin{align*}
I_1^2 &\leq C \int_{\calU} \bfrho(z)^{-2\gamma} \left( \int_{D(\bfi,R)^c} \frac{|f(z)-f(w)|^2 |g_k(w)|^2 \bfrho(w)^{2\gamma}}{|\bfrho(z,w)|^{n+1}} dV(w) \right) dV(z)\\
&\leq C \int_{D(\bfi,R)^c} |g_k(w)|^2 \left( \bfrho(w)^{2\gamma} \int_{\calU} \frac{|f(z)-f(w)|^2 \bfrho(z)^{-2\gamma}}{|\bfrho(z,w)|^{n+1}} dV(z) \right) dV(w)\\
&\leq C \varepsilon^2 \int_{D(\bfi,R)^c} |g_k(w)|^2 dV(w) \leq C \varepsilon^2 \|g_k\|^2.
\end{align*}
Let $M=\sup_k \|g_k\|$. Then we obtain
\[
I_1\leq  \sqrt{C} M \varepsilon .
\]
It remains to estimate $I_2$. By Minkowski's inequality, we have 
\begin{align*}
I_2 &\leq \int_{D(\bfi,R)} \left[\int_{\calU} |f(z)-f(w)|^2 |K(z,w)|^2 dV(z)\right]^{1/2} |g_k(w)| dV(w)\\
&= \int_{D(\bfi,R)} \|(f-f(w))k_w\| \sqrt{K(w,w)} |g_k(w)| dV(w).
\end{align*}
By Theorem \ref{thm:chaBO}, the function $\|(f-f(w))k_w\|$ is bounded on $\calU$. 
Also, it is obvious that the function $\sqrt{K(w,w)}$ is bounded on the compact set $\overline{D(\bfi,R)}$.
Therefore, there exists a positive constant $C_1$ depending on $R$ such that 
\[
I_2 \leq C_1 \sup_{w\in D(\bfi,R)} |g_k(w)|.
\]
Since $\{g_k\}$ converges to 0 uniformly on each compact subset of $\calU$, 
\[
I_2 \leq  \varepsilon
\]
for sufficiently large $k$.
Putting everything together, we conclude that $\|H_f g_k\|\to 0$ as $k\to\infty$ and hence $H_f$ is compact. 
\end{proof}

Combining Lemmas \ref{lem:VAHankel} and \ref{lem:VOHankel} with the decomposition
\[
VMO=VO+VA
\]
given in Theorem \ref{thm:BMO,VMO}, we see that if $f\in VMO$, then both $H_f$ and $H_{\overline{f}}$ are  compact. Thus, the necessity of Theorem \ref{thm:compactHankel} is proved. 


\section{The Bloch space}

In the last section, we identify the space of all holomorphic functions in $BMO$ with Bloch space on $\calU$. It should be pointed out that the definition of Bloch space of $\calU$ has never been clearly given in a direct way before, although the Siegel upper half-space and the unit ball are biholomorphically equivalent. To endow the Bloch space with property of M\"obius invariance, we shall first introduce the notion of invariant gradient. 

\subsection{The invariant gradient}
Let $\Omega \subset \mathbb{C}^n$ and $K$ be the Bergman kernel of $\Omega$. 
For $z\in\Omega$, define
\[
g_{i,j}(z) = \frac {1}{n+1} \frac {\partial^2 \log K(z,z)}{\partial z_i \partial \bar{z}_j}, \quad 1\leq i,j\leq n.
\]
Let $(g_{i,j})$ denote the $n\times n$ Hermitian matrix and $(g^{i,j})$ denote the inverse matrix. 
Set $g(z)=\det (g_{i,j}(z))$.
It is well known that the Laplace-Beltrami operator associated with the Bergman kernel $K$ is the differential operator $\widetilde{\Delta}$ defined by
\[
\widetilde{\Delta}  =\frac{2}{g} \sum_{i,j}\bigg\{\frac{\partial}{\partial \bar{z}_i}\bigg(g g^{i,j}\frac{\partial}{\partial z_j}\bigg)+\frac{\partial}{\partial z_j}\bigg(gg ^{i,j}\frac{\partial}{\partial \bar{z}_i}\bigg)\bigg\}.
\]
The operator $\widetilde{\Delta}$ is often referred to as the invariant Laplacian, since it has the property
\[
\widetilde{\Delta}(f\circ \psi)=(\widetilde{\Delta}f)\circ\psi \quad \text{for all} \quad \psi\in \mathrm{Aut}(\Omega),
\]
where $\mathrm{Aut}(\Omega)$ denotes the group of automorphisms of $\Omega$.
See \cite{Sto199} for instance. 
If $u$ and $v$ are $C^2$ functions, then we have
\[
\widetilde{\Delta}(uv)=u\widetilde{\Delta}v+2(\widetilde{\Delta}u)v+v\widetilde{\Delta}u,
\]
where $\widetilde{\nabla}u$ is the vector field defined by 
\[
\widetilde{\nabla}u=2\sum_{i,j} g^{i,j}\Bigg\{\frac{\partial u}{\partial \bar{z}_i} \frac{\partial}{\partial z_j}+
\frac{\partial u}{\partial z_j} \frac{\partial}{\partial \bar{z}_i} \Bigg\}.
\]
Again recall from \cite{Sto199} that 
\[
|\widetilde{\nabla}u|^2=(\widetilde{\nabla}u) \overline{u}=2\sum_{i,j} g^{i,j}\Bigg\{\frac{\partial u}{\partial \bar{z}_i} \overline{\frac{\partial u}{\partial \bar{z}_j}}+
\overline{\frac{\partial u}{\partial z_i}} \frac{\partial u}{\partial z_j}  \Bigg\}.
\]
In particular, if $f$ in $H(\Omega)$, the space of holomorphic functions in $\Omega$, then
\[
|\widetilde{\nabla}f|^2=2\sum_{i,j} g^{i,j} \overline{\frac{\partial f}{\partial z_i}} \frac{\partial f}{\partial z_j} .
\]
Furthermore, it's not hard to check that 
\[
\widetilde{\Delta} |f|^2=2 (\widetilde{\nabla} f)\overline{f}=2|\widetilde{\nabla} f|^2.
\]
Hence, 
\[
|\widetilde{\nabla}(f\circ\psi)|=|\widetilde{\nabla}f\circ\psi|\quad \text{for all} \quad \psi\in \mathrm{Aut}(\Omega)
\]
due to the invariance property of $\widetilde{\Delta}$. Therefore, the operator $|\widetilde{\nabla}|$ is usually called invariant gradient of $\Omega$.

\subsection{The definition of Bloch space}
With the above acknowledgements, we shall define the Bloch
space of $\calU$ in a direct way. Let $\Omega=\calU$. It is not hard to calculate that
\[
(g_{i,j}(z))= \frac {1}{\bfrho(z)^2} \left(
                               \begin{array}{cc}
                                 \bfrho(z) I_{n-1} + \bar{z}^{\prime} z^{\prime T}  & -(i/2) \bar{z}^{\prime} \\
                                 (i/2) z^{\prime T} & 1/4\\
                               \end{array}
                             \right)
\]
and 
\[
(g^{i,j}(z)) = \bfrho(z) \left(
                        \begin{array}{cc}
                          I_{n-1}  & 2i \bar{z}^{\prime} \\
                          -2i z^{\prime T} & 4\ImPt z_n\\
                        \end{array}
                       \right)
\]
for $z\in\calU$.
Then, on $\calU$, the Laplace-Beltrami operator is of the form
\[
4\bfrho(z) \bigg\{ \sum\limits_{j=1}^{n-1} \frac {\partial ^2}{\partial z_j \partial \bar{z}_j}
+ 2i \sum\limits_{j=1}^{n-1} \bar{z}_j \frac {\partial ^2}{\partial z_n \partial \bar{z}_j}
- 2i \sum\limits_{j=1}^{n-1} z_j \frac {\partial ^2}{\partial z_j \partial \bar{z}_n}
+ 4(\ImPt z_n) \frac {\partial ^2}{\partial z_n \partial \bar{z}_n} \bigg\}
\]
(see also \cite{Kor69}) and the invariant gradient has the form
\[
|\widetilde{\nabla}f(z)|^2=2\bfrho(z) \left( 4\bfrho(z) \left|\frac{\partial f}{\partial z_n}(z)\right|^2+\sum_{j=1}^{n-1} \left|\frac{\partial f}{\partial z_j}(z) +2i\bar{z}_j \frac{\partial f}{\partial z_n}(z)\right|^2 \right) 
\]
for $f\in H(\calU)$.

We now define the Bloch space of $\calU$, denoted by $\calB$, as the space of $f\in H(\calU)$ such that 
\[
\|f\|_{\calB}:=\sup\{|\widetilde{\nabla} f(z)|:z\in\calU\}<\infty.
\]
Clearly this is only a semi-norm and invariant under the action of $\text{Aut}(\calU)$. 

For the sake of the subject of this paper, we here are not going to make an intensive study of the Bloch space $\calB$ but only provide some properties to arrive at our main purpose of this section. In our another almost finished paper, more interesting properties of $\calB$ will be listed.

\subsection{Möbius transformations}
Let $\ball$ be the unit ball of $\mathbb{C}^n$.
It is known that the group $\mathrm{Aut}(\ball)$ is generated by the unitary transformations
on $\bbC^{n}$ along with the M\"obius transformations
$\varphi_{\xi}$ given by
\[
\varphi_{\xi}(\eta) := \frac {\xi-P_{\xi}\eta-(1-|\xi|^2)^{\frac {1}{2}}Q_{\xi}\eta}{1-\eta\cdot\overline{\xi}},
\]
where $\xi,\eta\in \ball$, $P_{\xi}$ is the orthogonal projection onto the space spanned
by $\xi$,
and $Q_{\xi}\eta=\eta- P_{\xi}\eta$.
It is easily shown that the mapping $\varphi_{\xi}$ satisfies
\[
\varphi_{\xi}(0)=\xi, \quad \varphi_{\xi}(\xi)=0, \quad \varphi_{\xi}(\varphi_{\xi}(\eta))=\eta.
\]
It is known that the invariant gradient of the unit ball is defined by
\[
|\widetilde{\nabla}_{\ball} h(\xi)|=|\nabla(h\circ\varphi_{\xi})(0)|
\]
for $h\in H(\ball)$ and $\xi\in\ball$, where
\[
\nabla h(\xi)= \left(\frac{\partial h}{\partial \xi_1}(\xi),\cdots,\frac{\partial h}{\partial \xi_n}(\xi) \right).
\]
See for instance \cite{Zhu05}.

Recall that $\calU$ is biholomorphically equivalent to $\ball$, via the Cayley transform $\Phi:\bbB\to \calU$ given by
\[
(z^{\prime}, z_{n})\; \longmapsto\; \left( \frac {z^{\prime}}{1+z_{n}},
i\left(\frac {1-z_{n}}{1+z_{n}}\right) \right)
\]
and its inverse
\[
\Phi^{-1}: \left(z^{\prime},z_{n}\right)\;\longmapsto\; \left(\frac {2iz^{\prime}}{i+z_{n}},
\frac {i-z_{n}}{i+z_{n}}\right).
\]
It is not hard to obtain their real Jacobians
\begin{equation*}\label{eqn:jacobian4phi}
\left(J_{R}\Phi\right)(\xi) = \frac {4}{|1+\xi_{n}|^{2(n+1)}}, \quad \xi\in\ball,
\end{equation*} 
and 
\begin{equation}\label{eqn:jacobian Phi^-1}
\left(J_{R}\Phi^{-1}\right)(z) = \frac {1}{4|\bfrho(z,\bfi)|^{2(n+1)}}, \quad z\in\calU.
\end{equation}
We refer to \cite[Chapter XII]{Ste93} for more properties of the Cayley transform. 

Through the Cayley transform, a class of M\"obius transformations of $\calU$ induced by $\varphi_{\xi}$ is given by
\[
\tau_z:=\Phi\circ \varphi_{\Phi^{-1}(z)}\circ\Phi^{-1}.
\]
Since obviously $\Phi(0)=\bfi$, $\tau_z(z)=\bfi$. 

\subsection{Estimates of invariant gradient}
By means of $\tau_z$, we obtain the relationship of invariant gradients of $\calU$ and $\ball$ as follows.
For any $f\in H(\calU)$, it is not hard to check that
\begin{equation}\label{eq:tildenablaf(i)}
|\widetilde{\nabla} f(\bfi)|=2|\nabla(f\circ\Phi)(0)|,
\end{equation}
then we have
\begin{align}
|\widetilde{\nabla}f(z)|&= |\widetilde{\nabla}(f\circ\tau_z)(\bfi)| = 2|\nabla(f\circ\tau_{z}\circ\Phi)(0)|\nonumber\\
&=2|\nabla(f\circ\Phi\circ\varphi_{\Phi^{-1}(z)})(0)|\nonumber\\
&=2|\widetilde{\nabla}_{\mathbb{B}}(f\circ\Phi)\big(\Phi^{-1}(z)\big)|\label{eq:gradient}
\end{align}
for all $z\in\calU$. 

\begin{lemma}\label{lem:nabla.Zhu 2.24a}
For $f\in H(\calU)$, then
\[
|f(z)-f(w)|\leq \sup_{u\in\gamma} |\widetilde{\nabla}f(u)| \beta(z,w)/2,
\]
where $\gamma$ is a geodesic joining $z$ to $w$ in the Bergman metric.
\end{lemma}
\begin{proof}
Let $\beta_{\ball}(\cdot,\cdot)$ denote the Bergman metric of $\ball$.
On one hand, \cite{Tim80I} that 
\[
|h(\xi)-h(\eta)| \leq  \sup_{\zeta\in\Gamma} |\widetilde{\nabla}_{\ball}h(\zeta)| \beta_{\ball}(\xi,\eta)
\]
for $h\in H(\ball)$ and $\xi,\eta\in\ball$, where $\Gamma$ is a geodesic joining $\xi$ to $\eta$ with respect to $\beta_{\ball}(\cdot,\cdot)$.
On the other hand, it is known that $\Phi$ sends the geodesic of $\ball$ with respect to $\beta_{\ball}(\cdot,\cdot)$ to the geodesic of $\calU$ with respect to $\beta(\cdot,\cdot)$, and 
\begin{equation}\label{eq:Bergmanmetric}
\beta(z,w)=\beta_{\ball}(\Phi^{-1}(z),\Phi^{-1}(w))
\end{equation}
for any $z,w\in\calU$; see \cite[Proposition 1.4.15]{Kra01} for instance. Therefore, let $z=\Phi(\xi)$, $w=\Phi(\eta)$, $u=\Phi(\zeta)$ and $\gamma=\Phi(\Gamma)$ be the geodesic joining $z$ to $w$ with respect to $\beta(\cdot,\cdot)$, we have
\begin{align*}
|f(z)-f(w)|&=| (f\circ\Phi)(\xi)-(f\circ\Phi)(\eta)|\\
&\leq \sup_{\zeta\in\Gamma} |\widetilde{\nabla}_{\ball}(f\circ\Phi)(\zeta)| \beta_{\ball}(\xi,\eta)\\
&=\sup_{u\in\gamma} \frac{|\widetilde{\nabla}f(u)|}{2} \beta(z,w),
\end{align*}
where the last equality uses \eqref{eq:gradient}, as desired.
\end{proof}

\begin{corollary}\label{cor:Blochbeta}
If $f\in\calB$, then
\[
|f(z)-f(w)|\leq \|f\|_{\calB} \beta(z,w)/2
\]
for all $z,w\in\calU$.
\end{corollary}

\begin{lemma}\label{lem:nabla.Zhu 2.24}
Suppose $r>0$ and $p>0$. There exists a positive constant $C$ such that 
\[
|\widetilde{\nabla}f(z)|^p \leq \frac{C}{|D(z,r)|}\int_{D(z,r)} |f(w)|^p dV(w)
\]
for all $f\in H(\calU)$ and all $z\in\calU$.
\end{lemma}
\begin{proof}
For $g\in H(\ball)$, there exists a positive constant $C$ depending on $r$ such that 
\[
|\nabla g(0)|^p \leq C \int_{D_{\ball}(0,r)} |g(\eta)|^p dV(\eta).
\]
See for instance \cite[Lemma 2.4]{Zhu05}. Then, for any $f\in H(\calU)$, we have 
\begin{align*}\label{eq:nabla.Zhu. 2.24}
&|\nabla(f\circ\Phi)(0)|^p \leq C  \int_{D_{\ball}(0,r)} |f\circ\Phi(\eta)|^p dV(\eta)\\
&\leq C \int_{D(\bfi,r)}  \frac{|f(w)|^p}{|\bfrho(w,\bfi)|^{2(n+1)}} dV(w)
\leq C \int_{D(\bfi,r)} |f(w)|^p  dV(w),
\end{align*}
where the second inequality uses \eqref{eqn:jacobian Phi^-1} and \eqref{eq:Bergmanmetric},  and the last inequality uses the fact that $|\bfrho(w,\bfi)|\geq 1/2$. 
This together with \eqref{eq:tildenablaf(i)} gives
\[
|\widetilde{\nabla} f(\bfi)|^p \leq C\int_{D(\bfi,r)} |f(w)|^p  dV(w).
\]
Notice that $\sigma_z^{-1}(D(\bfi,r))=D(z,r)$. Replacing $f$ by $f\circ\sigma_z^{-1}$ in the above inequality, we obtain
\[
|\widetilde{\nabla}f(z)|^p
\leq C \int_{D(\bfi,r)} |f\circ\sigma_z^{-1} (w)|^p  dV(w)
= \frac{C}{|D(z,r)|} \int_{D(z,r)} |f(w)|^p  dV(w),
\]
where the last equality uses \eqref{eqn:jacobian} and  \eqref{eqn:volBergmanball}. The proof of the lemma is complete.
\end{proof}

\subsection{Proof of Theorem \ref{thm:main3}}

\begin{theorem}\label{thm:BlochBMO}
Suppose $r>0$, $p\geq 1$ and $f\in H(\calU)$. Then the following conditions are equivalent:
\begin{enumerate}
\item[(a)] $f\in\calB$.
\item[(b)] There exists a positive constant $C$ such that
\[
\frac{1}{|D(z,r)|} \int_{D(z,r)} |f(w)-f(z)|^p dV(w) \leq C 
\]
for all $z\in\calU$.
\item[(c)] There exists a positive constant $C$ such that
\[
\frac{1}{|D(z,r)|} \int_{D(z,r)} |f(w)-\widehat{f}_r(z)|^p dV(w) \leq C 
\]
for all $z\in\calU$.
\item[(d)] There exists a positive constant $C$ with the property that for every $z\in\calU$ there is a complex number $c_z$ such that
\[
\frac{1}{|D(z,r)|} \int_{D(z,r)} |f(w)-c_z|^p dV(w) \leq C.
\]
\end{enumerate}
\end{theorem}
\begin{proof}
It follows from Lemma \ref{lem:nabla.Zhu 2.24a} that 
\begin{equation}\label{eq:thm:BlochBMO}
\frac{1}{|D(z,r)|}  \int_{D(z,r)}  |f(w)-f(z)|^p dV(w)\leq (r/2)^p \sup_{u\in D(z,r)} |\widetilde{\nabla}f(u)|^p.
\end{equation}
This proves that (a) implies (b).

Write
\[
f(w)-\widehat{f}_r(z)=f(w)-f(z)-(\widehat{f}_r(z)-f(z))
\]
and observe that
\[
\widehat{f}_r(z)-f(z)=\frac{1}{|D(z,r)|} \int_{D(z,r)} (f(w)-f(z)) dV(w).
\]
By H\"{o}lder's inequality, we have
\begin{align*}
|\widehat{f}_r(z)-f(z)|^p &\leq \frac{1}{|D(z,r)|^p}  \int_{D(z,r)} |f(w)-f(z)|^p dV(w) |D(z,r)|^{p(1-1/p)}\\
&=\frac{1}{|D(z,r)|} \int_{D(z,r)} |f(w)-f(z)|^p dV(w).
\end{align*}
Thus,
\begin{align}
&\frac{1}{|D(z,r)|} \int_{D(z,r)} |f(w)-\widehat{f}_r(z)|^p dV(w)\notag\\
&\leq 2^p\left( \frac{1}{|D(z,r)|} \int_{D(z,r)} |f(w)-f(z)|^p dV(w)+|\widehat{f}_r(z)-f(z)|^p\right)\notag\\
& \leq \frac{2^{p+1}}{|D(z,r)|} \int_{D(z,r)} |f(w)-f(z)|^p dV(w)\label{eq:thm:BlochBMOa}.
\end{align}
This proves that (b) implies (c).

That (c) implies (d) is trivial.

For any $z\in\calU$,  replacing $f$ by $f-c_z$ in Lemma \ref{lem:nabla.Zhu 2.24}, we obtain
\begin{equation}\label{eq:thm:BlochBMOb}
|\widetilde{\nabla}f(z)|^p \leq \frac{C}{|D(z,r)|} \int_{D(z,r)} |f(w)-c_z|^p dV(w).
\end{equation}
This proves that (d) implies (a).
The proof of the theorem is complete.
\end{proof}

\begin{theorem}\label{thm:Bloch_0VMO}
Suppose $r>0$, $p\geq 1$ and $f\in H(\calU)$. Then the following conditions are equivalent:
\begin{enumerate}
\item[(a)] $f\in\calB_0$.
\item[(b)] There is 
\[
\frac{1}{|D(z,r)|} \int_{D(z,r)} |f(w)-f(z)|^p dV(w) \to 0 \quad \text{as} \quad z\to\partial\widehat{\calU}.
\]
\item[(c)] There is
\[
\frac{1}{|D(z,r)|} \int_{D(z,r)} |f(w)-\widehat{f}_r(z)|^p dV(w) \to 0 \quad \text{as} \quad z\to\partial\widehat{\calU}.
\]
\item[(d)]  For every $z\in\calU$ there is a complex number $c_z$ such that
\[
\frac{1}{|D(z,r)|} \int_{D(z,r)} |f(w)-c_z|^p dV(w) \to 0 \quad \text{as} \quad z\to\partial\widehat{\calU}.
\]
\end{enumerate}
\end{theorem}
\begin{proof}
That (a) implies (b) follows from \eqref{eq:thm:BlochBMO}. That (b) implies (c) follows from \eqref{eq:thm:BlochBMOa}. That (c) implies (d) is trivial. That (d) implies (a) follows from \eqref{eq:thm:BlochBMOb}.
\end{proof}

For any $f\in H(\calU)$ and any $r>0$, by letting $p=2$ in Theorems \ref{thm:BlochBMO} and \ref{thm:Bloch_0VMO}, we see that $f\in\calB$ if and only if $f\in BMO_r$, $f\in\calB_0$ if and only if $f\in VMO_r$. 
Moreover, the quantities $\|f\|_{\calB}$ and $\|f\|_{BMO_r}$ are equivalent.
Therefore, Theorem \ref{thm:main3} is proved.


\begin{thebibliography}{99}

\bibitem{AFP88}
J. Arazy, D. Fisher and J. Peetre, Hankel operators on weighted Bergman spaces, Amer. J. Math. 110 (1988), 989-1053.

\bibitem{AFJP91}
J. Arazy, D. Fisher, S. Janson and J. Peetre, Membership of Hankel operators on the ball in unitary ideals, J. London Math. Soc. 43 (1991), 485-508.

\bibitem{Axl86}
S. Axler, The Bergman space, the Bloch space, and communicators of multiplication operators, Duke Math. J. 53 (1986), 315-332.


%



\bibitem{BCZ87}
C. Berger, L. Coburn and K. Zhu, BMO on Bergman spaces of the classical domains, Bull. Amer. Math. Soc. 17 (1987), 133-136.

\bibitem{BCZ88}
C. Berger, L. Coburn and K. Zhu, Function theory on Cartan domains and the Berezin-Toeplitz symbol calculus, Amer. J. Math. 110 (1988), 921-953.

\bibitem{BBCZ90}
D. B\'{e}koll\'{e}, C. Berger, L. Coburn and K. Zhu, BMO in the Bergman metric on bounded symmetric domains, J. Funct. Anal. 93 (1990), 310-350.



\bibitem{BY07}
H. Bommier-Hato and E. H. Youssfi, Hankel operators on weighted Fock spaces, Integral Equations Operator Theory 59 (2007), 1-17.






%
%

\bibitem{DK93}
M. M. Djrbashian and A. H. Karapetyan, Integral representations for some classes of functions holomorphic
in a Siegel domain, J. Math. Anal. Appl., 179 (1993), 91--109.



\bibitem{Gin64}
S. G. Gindikin, Analysis in homogeneous domains, Russian Math. Surveys 19(4) (1964), 1--89.

\bibitem{Har58}
P. Hartman, Completely continuous Hankel matrices, Proc. Amer. Math. Soc. 9 (1958), 862-866.



\bibitem{HW18}
Z. Hu and E. Wang, Hankel operators between Fock spaces, Integral Equations Operator Theory 90:37 (2018).

\bibitem{HL19}
Z. Hu and J. Lu, Hankel operators on Bergman spaces with regular weights, J. Geom. Anal. 29 (2019), 3494-3519.

\bibitem{JP92}
Q. Jiang and L. Peng, Toeplitz and Hankel type operators on the upper half-plane, Integral Equations Operator Theory 15 (1992), 744-767.


%



%
%
%

\bibitem{Kor69}
A. Kor\'{a}nyi, Harmonic functions on Hermitian hyperbolic space, Tran. Amer. Math. Soc. 135 (1969), 507-516.

\bibitem{Kra01}
S. Krantz, Function theory of several complex variables. Reprint of the 1992 edition.
AMS Chelsea Publishing, Providence, RI, 2001.



\bibitem{Li92}
H. Li, BMO, VMO and Hankel operators on the Bergman space of strongly pseudoconvex domains, J. Funct. Anal. 106 (1992), 375-408.

\bibitem{Li94}
H. Li, Hankel operators on the Bergman spaces of strongly pseudoconvex domains, Integral Equations Operator Theory 19 (1994), 458-476.


\bibitem{Liu2018}
C. Liu, Norm estimates for the Bergman and Cauchy-Szeg\"o projections over the Siegel upper half-space. Constr. Approx. 48 (2018), 385-413.


\bibitem{LLHZ18}
C. Liu, Y. Liu, P. Hu and L. Zhou, Two classes of integral operators over the Siegel upper half-space, Complex Anal. Oper. Theory (2019), no.3, 685-701.

\bibitem{LS20}
C. Liu and J. Si, Positive Toeplitz operators on the Bergman spaces of the Siegel upper half-space, Commun. Math. Stat. 8 (2020), 113-134.

%

%



%
%





\bibitem{Neh57}
Z. Nehari, On bounded bilinear forms, Ann. of Math. 65 (1957), 153-162.


\bibitem{PZZ16}
J. Pau, R. Zhao and K. Zhu, Weighted BMO and Hankel operators between Bergman spaces, Indiana Univ. Math. J. 65 (2016), 1639-1673.

\bibitem{Sch04}
G. Schneider, Hankel operators with antiholomorphic symbols on the Fock space, Proc. Amer. Math. Soc. 132 (2004), 2399-2409.

\bibitem{SY13}
K. Seip and E. H. Youssfi, Hankel operators on Fock spaces and related Bergman kernel estimates, J. Geom. Anal. 23 (2013), 170-201.

\bibitem{Ste93}
E. M. Stein, Harmonic analysis: real-variable methods, orthogonality, and oscillatory integrals.
Princeton University Press, Princeton, NJ, 1993.

\bibitem{Sto199}
M. Stoll, Invariant Potential Theory in the Unit ball of $\mathbb{C}^n$, London Mathematical Society Lecture Note Series 199, Cambridge University Press, 1994.

\bibitem{Tim80I}
R. Timoney,
Bloch functions in several complex variables I,
Bull. London Math. Soc. 12 (1980), 241-267.



\bibitem{WCZ15}
X. Wang, G. Cao and K. Zhu, BMO and Hankel operators on Fock-type spaces, J. Geom. Anal. 25 (2015), 1650-1665.




\bibitem{Zhu87}
K. Zhu, VMO, ESV, and Toeplitz operators on the Bergman space, Trans. Amer. Math. Soc. 302 (1987), 617-646.



\bibitem{Zhu05}
K. Zhu, Spaces of Holomorphic Functions in the Unit Ball, Graduate Texts in Math., vol. 226, Springer, New York 2005.




\end{thebibliography}
\end{document}